\documentclass[12pt,a4paper]{article}
\usepackage{amstext}
\usepackage{fancyhdr}
\usepackage{amsfonts,graphicx,bezier, amssymb}
\usepackage{amsmath}
\usepackage{caption}
\usepackage{mathtools}
\usepackage[utf8]{inputenc}
\usepackage{tikz}
\usepackage{tikz-cd}
\usetikzlibrary{shapes.geometric, arrows}
\usepackage{rotating}
\usetikzlibrary{arrows.meta}
\newenvironment{prf}{\noindent{\bf{Proof:}}~~}{\hfill\rule{1ex}{1ex}\vskip1.5ex}
\newcommand{\Z}{\mathbb Z}

\newcommand{\beqa}{\begin{eqnarray}}
\newcommand{\enqa}{\end{eqnarray}}
\newcommand{\beq}{\begin{eqnarray*}}
\newcommand{\enq}{\end{eqnarray*}}

\newtheorem{rem}{Remark}[section]

\newtheorem{cor}{Corollary}[section]
\newtheorem{propn}{Proposition}[section]
\newtheorem{defn}{Definition}[section]
\newtheorem{exam}{{\bf Example}}[section]
\newtheorem{thm}{Theorem}[section]

\newtheorem{lem}{Lemma}[section]

\newcommand{\noi}{\noindent}

\makeatletter
\providecommand*{\twoheadrightarrowfill@}{%
  \arrowfill@\relbar\relbar\twoheadrightarrow
}
\providecommand*{\twoheadleftarrowfill@}{%
  \arrowfill@\twoheadleftarrow\relbar\relbar
}
\providecommand*{\xtwoheadrightarrow}[2][]{%
  \ext@arrow 0579\twoheadrightarrowfill@{#1}{#2}%
}
\providecommand*{\xtwoheadleftarrow}[2][]{%
  \ext@arrow 5097\twoheadleftarrowfill@{#1}{#2}%
}
\makeatother

\begin{document}

\begin{center}
{\bf\Large Applications of reduced and coreduced modules II: Radicality of the  functor $\text{Hom}_R(R/I, -)$}
\end{center}

\vspace*{0.3cm}
\begin{center}
David Ssevviiri\\
\vspace*{0.3cm}
Department of Mathematics\\
Makerere University, P.O BOX 7062, Kampala, Uganda\\
Email: david.ssevviiri@mak.ac.ug, ssevviiridaudi@gmail.com
\end{center}

\begin{abstract}This is the second in a series of papers  highlighting  the applications of reduced and coreduced modules. Let $R$ be a commutative unital ring and $I$ be an ideal of $R$.  We give  necessary and sufficient conditions in terms of   $I$-reduced $R$-modules   for the functor $\text{Hom}_R(R/I, -)$ on an Abelian full subcategory of the category of $R$-modules to be a  radical. $I$-reduced  and $I$-coreduced $R$-modules   provide  a natural setting  for a generalisation of Jans' correspondence, and lead to the construction of a new radical class of rings.
\end{abstract}


\vspace*{0.4cm}

{\bf MSC 2010} Mathematics Subject Classification:  16S90, 13D30,  16D80, 13D07

 \section{Introduction}

 \paragraph{\bf Why should we care about radicals?}  Throughout the paper, $R$ is  a commutative unital ring and $I$ is an ideal of $R$. Radicals play an important role in describing the structure of rings and modules. For a comprehensive treatment of radicals of rings, see \cite{Wiegandt}; and for modules, see  \cite{Stenstrom}. The most interesting and useful radicals are those which are left exact. If a radical (which by definition is a functor) on the category $R$-Mod of all $R$-modules  is  left exact, then it is  in a one-to-one correspondence with a hereditary torsion theory which in turn is in a one-to-one correspondence with a Gabriel topology, \cite[Chapter VI, Theorem 5.1]{Stenstrom}. To every left exact radical $\gamma$ on $R$-Mod and to every $M\in R$-Mod, we can  associate a module of fractions of $M$, given by $$\mathcal{M}_{\gamma}:=\lim_{\stackrel{\rightarrow}{J\in \mathcal{G}}}\text{Hom}_R(J, M/\gamma(M)),$$  where $\mathcal{G}$ is the Gabriel topology associated  with $\gamma$, \cite[page 197]{Stenstrom}. This is a general framework for constructing modules of quotients. It turns out that the usual local modules at a prime ideal are just a special case.   In general, torsion theories can be studied for any Abelian category, \cite{Dickson}. In the setting of triangulated categories, the analogue of torsion theories is what is called  {\it $t$-structures} \cite{BBD}; a common theme in  both  algebra and algebraic geometry.

 \paragraph{\bf Reduced modules and radicals in the literature.}\noi There is an intimate relationship between rings and modules which are reduced and some of their radicals.  A ring is reduced if and only if, its prime radical (also called the nilradical) is zero. An $R$-module $M$ is reduced if and only if, the locally nilradical $a\Gamma_a(M)$ (introduced in \cite{Kyom}) vanishes for all $a\in R$, where $\Gamma_a(M)$ is  the $(a)$-torsion submodule of $M$ and $a\Gamma_a(M)$ is the left multiplication by $a$ of the submodule $\Gamma_a(M)$, see \cite[Proposition 2.3]{Kyom}.

 \paragraph{\bf Jans' correspondence and how it is genralised in the paper.}
 In  \cite{Jan}, Jans showed that the map $I\mapsto T_I$ defines a one-to-one correspondence between an idempotent ideal $I$ of $R$ and the torsion-torsionfree (TTF) class $T_I$ given by $T_I:=\{M\in R\text{-Mod}~:~IM=0\}$.  In Theorem \ref{main}, we generalise this correspondence.  For an ideal $I$ of $R$ which is flat as an $R$-module, if  $\mathcal{A}_I$  and $\mathcal{B}_I$ denote Serre subcategories of $R$-Mod consisting of $I$-reduced $R$-modules and $I$-coreduced $R$-modules respectively, then we define a torsion theory $(\mathcal{T}, \mathcal{F}) $ on $\mathcal{A}_I$ and a torsion theory $(\mathfrak{T}, \mathcal{T}) $ on $\mathcal{B}_I$ such that $\mathcal{T}$ is simultaneously a torsion class in $(\mathcal{T}, \mathcal{F}) $ and a torsionfree class in $(\mathfrak{T}, \mathcal{T}) $. It turns out that Jans' correspondence is a special case, i.e., whenever $I^2=I$, we have $\mathcal{A}_I=\mathcal{B}_I=R$-Mod and the two torsion theories coincide with those given by Jans.

 \paragraph{\bf The conceptual contribution.}
  The main conceptual contribution of the paper is that the failure of Jans' correspondence for non-idempotent ideals is exactly the failure of $Hom_R(R/I, -)$ to be a radical on $R$-Mod. Restricting to suitable categories of $I$-reduced modules restores radicality and yields generalised torsion-theoretic correspondences.

 \paragraph{\bf Some consequences:  a radical class and applications to spectral sequences}
  Interpreting the module-theoretic results above in the setting of rings leads to a radical class (the analogue of a torsion class for modules). Let $R$ be a fixed commutative ring and $I$ an ideal of $R$. We consider the category $\mathcal{C}_R$  whose objects are rings equipped with a compatible $R$-module structure and whose morphisms are ring homomorphisms preserving the $R$-module structure.  We show in Theorem \ref{rclass} that for an  idempotent ideal $I$ of a ring $R$, the class of rings $\Psi_I:=\{S\in \mathcal{C}_R~:~ IS=0\}$ is a radical class. Furthermore, for any idempotent ideal $I$ of $R$, $\text{Hom}_R(R/I,R)\cong \Psi_I(R)=\Gamma_I(R)$, see Corollary \ref{crad}.
 Lastly, we utilise properties of $I$-reduced and $I$-coreduced modules to compute the Grothendieck spectral sequences associated with the $I$-torsion functor $\Gamma_I$ and  the $I$-adic completion functor $\Lambda_I$.
  Note that other applications of reduced and coreduced modules have already appeared in \cite{Kimuli} and \cite{David}. In \cite{Kimuli}, reduced modules were used to characterise regular modules. In \cite{David}, we demonstrated that $I$-reduced and $I$-coreduced modules provide a setting for which the Greenlees-May Duality and the Matlis-Greenlees-May Equivalence hold in the category of $R$-modules. Extension of the same to the category of chain complexes of $R$-modules and to the derived category of $R$-Mod was done in \cite{Amanuel}.

 \paragraph{\bf Organisation of the paper.} The paper  is organised  into six sections.   In section 2, we lay out  the   tools  required in the proofs of the main results.    In section 3, one of  the main results, namely Theorem \ref{main}, is proved.
  In section 4, we refine some results that already exist in the literature. The radical class of rings arising from the torsion theories studied is given in section 5, see Theorem \ref{rclass}.  In the last section, section 6, we give the spectral sequences associated with the functors $\Gamma_I$ and $\Lambda_I$ in the context of $I$-reduced and $I$-coreduced modules.

\section{Reduced modules $\&$ the functor $\text{Hom}_R(R/I, -)$}
\begin{defn} \label{red}\rm\cite{David}
 Let $I$ be an ideal of a ring $R$. An $R$-module $M$ is

 \begin{enumerate}
  \item
 {\it $I$-reduced} if for all   $m\in M$, $I^2m=0$ implies that $Im=0$;
 \item {\it $I$-coreduced} if $I^2M= IM$.
 \end{enumerate}

 \end{defn}

\paragraph{\bf The chronology and related notions.}
An $R$-module is {\it reduced} (resp. {\it coreduced}) if it is $I$-reduced (resp. $I$-coreduced) for all ideals $I$ of $R$.  The $R$-module $R$ is reduced if and only if, the ring $R$ is reduced. Reduced modules were introduced by Lee and Zhou in \cite{Lee}. Coreduced modules were first defined by   Ansari-Toroghy and Farshadifar in \cite{Ansari}  where they were called semisecond modules.  Generalised versions of $I$-reduced  $R$-modules exist in the literature  with different names. They are called;  modules  with bounded $I$-torsion, generalised $I$-reduced modules and
modules whose sequence of submodules $\{(0:_MI^k)\}_{k\in\Z^+}$ is stationary. See for instance \cite[Definition 5.5]{yk} $\&$ \cite[Sec. 7]{Rohrer}, \cite{Kyom3} and \cite[Example 7.3.2 (c), Proposition 3.1.10]{Peter} respectively.  On the other hand, $I$-coreduced modules were introduced in \cite[Theorem 2.3]{Nam} where they are called  modules for which the chain  $\{I^tM\}_{t\in \Z^+}$ of submodules of $M$ is stationary.

\begin{defn}\rm
A functor $\gamma: R\text{-Mod}\rightarrow R\text{-Mod}$ which associates to every $R$-module $M$, a submodule $\gamma(M)$ of $M$ is a:
\begin{itemize}
 \item[(i)] {\it preradical} if  for every $R$-homomorphism $f: M\rightarrow N$, $f(\gamma(M))\subseteq \gamma(N)$;
 \item[(ii)]   ~{\it radical} if  it is a preradical and for all $M\in R$-Mod, $\gamma(M/\gamma(M))=0$.
\end{itemize}
\end{defn}

\begin{exam}\rm The following are some examples of radicals defined on the category of $R$-modules.
 \begin{enumerate}
  \item For any $R$-module $M$, the intersection of all maximal submodules of $M$ is an idempotent radical called the {\it Jacobson radical} of $M$.

  \item Let $I$ be an ideal of $R$. The functor $\delta_I : R\text{-Mod}\rightarrow  R\text{-Mod}$ which associates to every $R$-module $M$, the submodule $IM$, is a radical.

  \item   For any finitely generated ideal $I$ of a ring $R$, the $I$-torsion functor $\Gamma_I$  is a left exact idempotent radical on $R$-Mod. It associates to every $R$-module $M$, the submodule $\Gamma_I(M):=\{m\in M~:~ I^km=0~\text{for some} ~k\in \Z^+\}$ of $M$.

  \item Let $S$ be a multiplicatively closed subset of an integral domain $R$. The submodule $$t(M):=\{m\in M~:~sm=0~\text{for some}~s\in S\}$$ of $M$ defines a left exact idempotent radical of $R$-Mod.

  \item For any $R$-module $M$,  the Bass torsion $B(M):= \text{Ker}(M\rightarrow M^{**})$, where $M^*:=\text{Hom}_R(M, R)$, is a radical. This radical first appeared in \cite{Bass} and was recently studied in \cite{Alex1} and \cite{Alex2}.

  \item To every $R$-module $S$, there is a left exact radical $r_S$ which associates to every $R$-module $M$,  the submodule  $$r_S(M):=\{m\in M~:~f(m)=0~\text{for all}~f:M\rightarrow E\},$$ where $E$ is the injective hull of $S$, see \cite{Sim}.
 \end{enumerate}

\end{exam}

\paragraph\noi
 Left exact radicals are also called {\it idempotent kernel functors} in some literature, see for instance, \cite{Sim}.   Denote  the statement ``$N$ is a submodule of $M$'' by $N\leq M$ and for any ideal $I$ of a ring $R$, denote the  submodule $\{m\in M~:~Im=0\}$ of $M$ by $(0:_MI)$.

\begin{propn}\label{prd}
 For any ideal $I$  of a ring $R$, the functor $\text{Hom}_R(R/I, -)$ on  the category $R$-Mod is a preradical.
\end{propn}

\begin{prf}
 Let $M$ and $N$ be $R$-modules and $f: M\rightarrow N$ be an $R$-homomorphism. Since $\text{Hom}_R(R/I, M)$ is naturally isomorphic to  $(0:_MI)$ for any $M\in R$-Mod, it is enough to show that $f((0:_MI))\subseteq (0:_NI)$. Let $m\in (0:_MI)$, then $Im=0$ and $If(m)= f(Im)=0$. This shows that $f(m)\in (0:_NI)$. Hence, $f((0:_MI))\subseteq (0:_NI)$.
\end{prf}

\begin{lem}\label{L2.1}
 Let $M$ be an $R$-module and $I$ be an  ideal of $R$. For any positive integer $k$,
 $$\left(0:_{\frac{M}{(0:_MI^k)}}I\right)\cong \frac{(0:_MI^{k+1})}{(0:_MI^k)}.$$
\end{lem}

\begin{prf}Define  a map $f: (0:_MI^{k+1})\rightarrow \left(0:_{\frac{M}{(0:_MI^k)}}I\right)$ by $m\mapsto m + (0:_MI^k)$. $f$ is an $R$-epimorphism with  kernel  $(0:_MI^k)$. The desired result becomes immediate by applying the first isomorphism theorem.
\end{prf}

\paragraph{\bf A module via apolarity.}
 Let $k$ be an algebraically closed field of characteristic zero. If $R:=k[x_1,\cdots, x_n]$ and $S:=k[X_1, \cdots, X_n]$, then $S$ is an $R$-module via the following action, called {\it apolarity} or {\it contraction}. $$R\circ S\rightarrow S$$

 $$(x^{\alpha}, X^{\beta})\mapsto  x^{\alpha}\circ X^{\beta}:=\begin{cases} \frac{\beta !}{(\beta -\alpha)!}X^{\beta-\alpha}, ~\text{if}~\beta_i\geq \alpha_i ~\text{for all}~i, \cdots, n\cr
 0, ~\text{otherwise}
 \end{cases};$$ where $x^{\alpha}= x_1^{{\alpha}_1}, \cdots, x_n^{{\alpha}_n}$ and $X^{\beta}= X_1^{{\beta}_1}, \cdots, X_n^{{\beta}_n}$. Usually, the apolarity action  shows up whenever one is working with Macaulay inverse systems,  see \cite{Elias, Iarrobino}  among others.

 \begin{exam}\rm \label{exapolarity}
  Let $M:=k[X]= k\oplus kX\oplus kX^2\oplus \cdots$ and $R:=k[x]$. $M$ is an $R$-module under apolarity action. If $I:=\langle x^2 \rangle$ is the ideal of $R$ generated by $x^2$, then $(0:_MI)= k\oplus kX$ and  $(0:_MI^2)= k\oplus kX\oplus kX^2\oplus kX^3$ so that, $(0:_MI^2)/ (0:_MI)= kX^2\oplus kX^3~ \text{mod}~ (k\oplus kX)$. Furthermore, $M/(0:_MI)= kX^2\oplus kX^3\oplus kX^4\oplus \cdots~\text{mod}~(k\oplus kX)$ and $(0:_{M/(0:_MI)}I)=  kX^2\oplus kX^3 ~\text{mod}~(k\oplus kX)$. It is then clear that  $(0:_{M/(0:_MI)}I)\cong (0:_MI^2)/ (0:_MI)$.
 \end{exam}

\begin{lem}\label{lr} Let $I$ be an ideal of a ring $R$. If  the functor $\Gamma_I$ is a radical,  then  it is the smallest radical containing (a naturally isomorphic copy of) the preradical $\text{Hom}_R(R/I, -)$.
\end{lem}

\begin{prf}
 Let $$\gamma(M):= (0:_M I),$$ and $\gamma_2(M)$ be the  unique submodule of $M$  containing $\gamma(M)$ such that $$\frac{\gamma_2(M)}{\gamma(M)}=\gamma\left(\frac{M}{\gamma(M)}\right).$$ It follows by Lemma \ref{L2.1} that $$\gamma\left(\frac{M}{\gamma(M)}\right)=\left(0:_{\frac{M}{(0:_MI)}}I\right)\cong\frac{(0:_MI^2)}{(0:_MI)}$$ and $$\gamma_2(M)= (0:_MI^2).$$

 If $k$ is not a limit ordinal, define $\gamma_{k}(M)$ to be the unique submodule of $M$  containing $\gamma_{k-1}(M)$ such that $$\frac{\gamma_{k}(M)}{\gamma_{k-1}(M)}=\gamma\left(\frac{M}{\gamma_{k-1}(M)}\right).$$

 $$\gamma\left(\frac{M}{\gamma_{k-1}(M)}\right)=\left(0:_{\frac{M}{\left(0:_MI^{k-1}\right)}}I\right)\cong \frac{(0:_MI^{k})}{(0:_MI^{k-1})}.$$ So, $$\gamma_k(M)= (0:_M I^k).$$

 If $k$ is a limit ordinal, then define $\gamma_{k}(M)$ to be the submodule of $M$ given by
 $$ \gamma_{k}(M):=\sum_{r<k}\gamma_r(M).$$

 It is well known that if $\gamma$ is a preradical, then the process described above leads to an ascending chain of submodules

 \begin{equation} \label{chain}
  \gamma(M)\subseteq \gamma_2(M)\subseteq \cdots \subseteq \gamma_k(M)\subseteq \cdots
 \end{equation}

 which terminates at a radical and this radical is the smallest containing the preradical $\gamma$, see \cite[Chapter VI, Proposition 1.5]{Stenstrom}. As seen above, Chain (\ref{chain}) is nothing but the chain

 \begin{equation}\label{chain2}
  (0:_MI)\subseteq (0:_MI^2)\subseteq \cdots \subseteq (0:_MI^k) \subseteq \cdots
 \end{equation}

  whose union is the $R$-module $\Gamma_I(M)$.
 So, if
 $\Gamma_I(-)$ is a radical, then it is the smallest radical containing $(0:_-I)\cong\text{Hom}_R(R/I, -)$.
\end{prf}

\begin{rem} \rm\label{R1}
 Whereas Stenstrom in \cite{Stenstrom} defines radicals and torsion theories mainly on the category $R$-Mod, the two notions also make sense in the setting of any Abelian category, \cite{Dickson}.
\end{rem}

\begin{rem}\rm
 The obstruction to extending Jans' correspondence to non-idempotent ideals is that $\text{Hom}_R(R/I,-)$ fails to be a radical on $R$-Mod. Restricting to $I$-reduced modules restores radicality.
\end{rem}

\paragraph{\bf Serre subcategories.}
Let \(\mathcal A\) be an abelian category. A full subcategory $\mathcal S\subseteq \mathcal A$ is called a \emph{Serre subcategory} if for every short exact sequence $0\longrightarrow A' \longrightarrow A
\longrightarrow A'' \longrightarrow 0$ in \(\mathcal A\), one has $A\in \mathcal S$ if and only if, $A',A''\in \mathcal S$.
Equivalently, a full subcategory \(\mathcal S\) is a Serre subcategory if it is closed under subobjects, quotient objects, and extensions. Since every object in an abelian category is built from subobjects, quotients, and extensions, a Serre subcategory is stable under all natural operations arising from exact sequences.

\begin{propn}\label{t3}
For any ideal $I$   of a ring $R$, let $\mathcal{A}_I$ and $\mathcal{B}_I$ be Serre full subcategories of $R$-Mod.

\begin{enumerate}
\item The functor  $\text{Hom}_R(R/I, -)$ on $\mathcal{A}_I$ is a  radical if and only if, $\mathcal{A}_I$ consists of $I$-reduced $R$-modules.

\item     The radical $\delta_I$ which associates to every $R$-module $M$ in $\mathcal{B}_I$, the $R$-submodule $\delta_I(M):=IM$ is  idempotent if and only if, $\mathcal{B}_I$ consists of $I$-coreduced $R$-modules.
\end{enumerate}

\end{propn}

\begin{prf}

\begin{enumerate}
 \item  In light of Lemma \ref{lr} and Remark \ref{R1}, the preradical $\text{Hom}_R(R/I, -)$ is a radical on $\mathcal{A}_I$ if and only if, for all $M\in \mathcal{A}_I$, $\text{Hom}_R(R/I, M)\cong (0:_MI)\subseteq (0:_MI^2)\subseteq (0:_MI^3)\subseteq \cdots$ is a constant chain of submodules of $M$ if and only if, for all $M\in \mathcal{A}_I$, $M$ is an $I$-reduced $R$-module. $\mathcal{A}_I$ being a Serre subcategory ensures that the $I$-reduced $R$-module $\text{Hom}_R(R/I, M)\cong (0:_MI)$ which is $I$-reduced (since it is annihilated by $I$) belongs to $\mathcal{A}_I$. Note that every radical is an endofunctor.

 \item The radical $\delta_I: \mathcal{B}_I\rightarrow \mathcal{B}_I$, $M\mapsto \delta_I(M):=IM$ is idempotent if and only if, for all $M\in \mathcal{B}_I$, $I^2M=IM$ if and only if, for all $M\in  \mathcal{B}_I$, $M$ is an $I$-coreduced $R$-module. Just like in the proof of part 1, whereas in this special case the submodule $IM$ is $I$-coreduced, the requirement that the subcategory $\mathcal{B}_I$ be  Serre ensures that $IM$ belongs to $\mathcal{B}_I$.

\end{enumerate}
\end{prf}

\begin{exam}\rm
 The preradical $\text{Hom}_R(R/I, -)$ on $R$-Mod need not be a radical. In  Example \ref{exapolarity}, we have an $R$-module $M$ and an ideal $I$ of $R$ such that the inclusion $(0:_MI)\subseteq (0:_MI^2)$ is strict. There exists an element $m\in M$ such that $I^2m=0$ but $Im\not=0$. Thus, $M$ is not $I$-reduced and by Proposition \ref{t3}, $\text{Hom}_R(R/I, -)$ is not a radical on $R$-Mod.
 \end{exam}

 \paragraph{\bf Effect on the $I$-torsion and  $I$-adic completion functors.}
  The  functors $
\Gamma_I ~:~R\text{-Mod}\rightarrow R\text{-Mod}$, $$M\mapsto \Gamma_I(M):= \lim_{\stackrel{\rightarrow}{k}}\text{Hom}_R(R/I^k, M)$$ and $ \Lambda_I ~:~R\text{-Mod}\rightarrow R\text{-Mod}$,
$$M\mapsto \Lambda_I(M):=\lim_{\stackrel{\leftarrow}{k}} (M/I^kM)$$ are called the {\it $I$-torsion functor} and the {\it $I$-adic completion functor} respectively.   An $R$-module is {\it $I$-torsion} (resp. {\it $I$-complete}) if $\Gamma_I(M)\cong M$ (resp., $\Lambda_I(M)\cong M$). Let $(R\text{-Mod})_{I\text{-red}}$ (resp. $(R\text{-Mod})_{I\text{-cor}}$) denote the subcategory of all $I$-reduced (resp. $I$-coreduced) $R$-modules.
 If $M\in (R\text{-Mod})_{I\text{-red}}$, then  $\Gamma_I(M)\in (R\text{-Mod})_{I\text{-cor}}$, and if $N\in (R\text{-Mod})_{I\text{-cor}}$, then  $\Lambda_I(N)\in (R\text{-Mod})_{I\text{-red}}$. By \cite[Theorem 1.1]{David}, the functors $$\Gamma_I : (R\text{-Mod})_{I\text{-red}}\rightarrow (R\text{-Mod})_{I\text{-cor}}$$  and $$\Lambda_I : (R\text{-Mod})_{I\text{-cor}}\rightarrow (R\text{-Mod})_{I\text{-red}}$$ form an adjoint pair.  A functor $\Omega : R\text{-Mod}\rightarrow R\text{-Mod}$ is a {\it coradical} if
 there exists a radical $\gamma: R\text{-Mod}\rightarrow R\text{-Mod}$ such that for every $M\in R\text{-Mod}$, $\Omega(M)= M/\gamma(M)$. We restrict this definition to the setting of our Serre subcategories $\mathcal{A}_I$ and $\mathcal{B}_I$.

\begin{propn} Let $\mathcal{A}_I$ and $\mathcal{B}_I$ be Serre full subcategories of $R$-Mod which  consist of $I$-reduced $R$-modules and $I$-coreduced $R$-modules respectively.

 \begin{enumerate}
  \item The functor $\Gamma_I: \mathcal{A}_I\rightarrow \mathcal{A}_I$ is a left exact radical.

  \item The functor $\Lambda_I: \mathcal{B}_I\rightarrow \mathcal{B}_I$ is a right exact coradical.
 \end{enumerate}

\end{propn}

\begin{prf}
 \begin{enumerate}
  \item Since the functor $\Gamma_I$  restricted to $\mathcal{A}_I$ is naturally isomorphic to the left exact functor $\text{Hom}_R(R/I,-)$, see \cite[Proposition 2.2(4)]{David}. By
Proposition \ref{t3}, the functor   $\Gamma_I: \mathcal{A}_I\rightarrow \mathcal{A}_I$  is a left exact radical.

  \item Since for $M\in \mathcal{B}_I$, $\Lambda_I(M)\cong M/IM$, \cite[Proposition 2.3(4)]{David} and by Proposition \ref{t3} the assignment $M\mapsto IM$ is a radical, the functor $\Lambda_I: \mathcal{B}_I\rightarrow \mathcal{B}_I$ is a coradical. It is right exact because $\Lambda_I$ on $\mathcal{B}_I$ is isomorphic to $R/I\otimes -$ which is right exact.
 \end{enumerate}

\end{prf}

\section{Torsion theories induced by a flat ideal}

\paragraph{}
 A {\it torsion theory} $\tau$ for an Abelian category $\mathcal{C}$ is a pair $(\mathcal{T}, \mathcal{F})$ of classes of objects of  $\mathcal{C}$ such that
 \begin{enumerate}
  \item $\text{Hom}(T, F)=0$ for all $T\in \mathcal{T}$, $F\in \mathcal{F}$;
  \item if $\text{Hom}(A, F)=0$ for all  $F\in \mathcal{F}$, then $A\in \mathcal{T}$;

  \item if $\text{Hom}(T, B)=0$ for all $T\in \mathcal{T}$,  then $B \in \mathcal{F}$.
  \end{enumerate}

  $\mathcal{T}$ is called the {\it torsion class} of $\tau$ and its objects are called {\it torsion objects}, whereas $\mathcal{F}$ is called the {\it torsionfree class} of $\tau$ and its objects are called {\it torsionfree objects}. A class $\mathcal{H}$ of objects of an Abelian category $\mathcal{C}$ is a {\it torsion-torsionfree} class (written TTF class for brevity) if it is  simultaneously a torsion class for one torsion theory and a torsionfree class for another. In other words, there exist torsion theories $(\mathcal{T}_1,\mathcal{F}_1)$ and $(\mathcal{T}_2,\mathcal{F}_2)$ in $\mathcal{C}$ such that  $\mathcal{H}=\mathcal{T}_1 = \mathcal{F}_2$.  A torsion class is {\it hereditary} if it is closed under taking submodules.

\begin{thm}{\rm\cite[Corollary 2.2]{Jan}[{\bf Jans' Correspondence}]}
Let $R$ be a ring.  There is a one-to-one correspondence between  idempotent ideals $I$ of $R$ and  TTF classes within $R$-Mod given by $I\mapsto \{M\in R\text{-Mod}~: IM=0\}$.
\end{thm}


\paragraph{}\noi
Consider the full subcategories of $R$-Mod below:
$$\mathcal{T}:= \left\{ M\in R\text{-Mod}~:~(0:_MI)=M \right\}$$   and   $$\mathcal{F}:= \left\{ M\in R\text{-Mod}~:~(0:_MI)=0 \right\}.$$ $\mathcal{T}$ is not closed under extensions and the functor $\text{Hom}_R(R/I, -)$ on $R$-Mod which is naturally isomorphic to $(0:_{-}I)$ is a left exact idempotent preradical which is  not a radical in general. Accordingly, $\mathcal{T}$ is just a pretorsion class, see \cite[Chapter VI, Proposition 1.4]{Stenstrom}.

\paragraph{} If $I$  is an idempotent ideal of $R$, then the functor $\text{Hom}_R(R/I, -)$ on $R$-Mod becomes a radical. In addition, $\mathcal{T}$ is closed under extensions and therefore  becomes a TTF class. Conversely, if $\mathcal{T}$ is a TTF class in $R$-Mod, then $I$ is an idempotent ideal of $R$. So, we get a one-to-one correspondence between idempotent ideals $I$ of $R$ and TTF classes within $R$-Mod given by $I\mapsto \left\{ M\in R\text{-Mod}~:~(0:_MI)=M \right\}$, see \cite[Corollary 2.2]{Jan}. We call this Jans' correspondence. Theorem \ref{main} below demonstrates that if the condition of idempotence of the ideal $I$ is relaxed to $I$ being a flat $R$-module, then we obtain still two torsion theories but on Serre subcategories  $\mathcal{A}_I$ and $\mathcal{B}_I$ respectively.

\begin{thm}\label{main}{\rm{\bf  [Generalisation of Jans' Correspondence]} }Let $I$ be an ideal of a ring $R$ such that $I$ is flat as an $R$-module.
Let $\mathcal A_I$ and $\mathcal B_I$   be Serre subcategories of $R$-Mod consisting of $I$-reduced and $I$-coreduced  $R$-modules respectively. Define
$$\mathcal{T}:= \left\{ M\in \mathcal{A}_I~:~(0:_MI)=M \right\}, $$      $$\mathcal{F}:= \left\{ M\in \mathcal{A}_I~:~(0:_MI)=0 \right\}$$ and
$$\mathfrak{T}:= \left\{ M\in \mathcal{B}_I~:~IM=M \right\}. $$ Then the following torsion theories arise:

\begin{enumerate}

 \item The pair  $(\mathcal{T}, \mathcal{F})$   forms a  hereditary  torsion theory  in the Abelian subcategory $\mathcal{A}_I$.

 \item The pair    $(\mathfrak{T}, \mathcal{T})$    forms a  torsion theory  in the Abelian  subcategory   $\mathcal{B}_I$.

 \item Consequently,  $\mathcal T$ occurs simultaneously as the torsion class in the torsion theory $(\mathcal T,\mathcal F)$ on $\mathcal A_I$ and as the torsionfree class in the torsion theory
$(\mathfrak T,\mathcal T)$ on $\mathcal B_I$.
\end{enumerate}

\end{thm}

\begin{prf}\rm
 \begin{enumerate}
  \item It is easy to see that $\mathcal{T}$ is closed under submodules, quotients, direct sums and under taking direct products in $\mathcal{A}_I$. We only prove that it is closed under  extensions. Let $A, B$ and $C$ be $R$-modules and $0\rightarrow A \rightarrow B \rightarrow C\rightarrow 0$ be a short exact sequence. Since $I$ is flat, $IM\cong I \otimes M$ for any $R$-module $M$ and the functor $I\otimes_R-$ is exact. Hence from $ 0\to A\to B\to C\to0$, we obtain an exact sequence $0\to IA\to IB\to IC\to0$. If $A,C\in\mathcal T$, then $IA=0$ and $IC=0$. Hence
$0\to0\to IB\to0\to0$ is exact, so $IB=0$ and  $B\in\mathcal T$.
This shows that $\mathcal T$ has the closure properties expected of a hereditary torsion class and of a torsionfree class. The corresponding torsion theories are constructed in parts 2 and 3 below.

\item For any $R$-module $M$, the submodule $(0:_MI)\cong \text{Hom}_R(R/I, M)$ of $M$  is annihilated by $I$. Hence, it is an $I$-reduced $R$-module.  By the fact that $\mathcal{A}_I$ is Serre, the submodule  $(0:_MI)$ belongs to $\mathcal{A}_I$.  Restricting to the subcategory $\mathcal{A}_I$, we get a functor $\text{Hom}_R(R/I, -):  \mathcal{A}_I\rightarrow \mathcal{A}_I$ which by Proposition \ref{t3} is a left exact radical. By \cite[Chapter VI, Proposition 3.1]{Stenstrom},  this radical is in a one-to-one correspondence with the hereditary torsion theory $(\mathcal{T}, \mathcal{F})$, where
 $ \mathcal{T}= \left\{ M\in \mathcal{A}_I~:~(0:_MI)=M \right\}$ and        $\mathcal{F}= \left\{ M\in \mathcal{A}_I~:~(0:_MI)=0 \right\}$.

 \item First observe that an $R$-module $M$ with the property that $IM=0$  for an ideal $I$ of $R$ is both $I$-reduced and $I$-coreduced and that $IM=0$ if and only if, $(0:_MI)=M$. This implies that the following three classes of $R$-modules coincide: $\left\{ M\in \mathcal{A}_I~:~IM=0 \right\}$, $\left\{ M\in \mathcal{B}_I~:~IM=0 \right\}$ and  $\left\{ M\in R\text{-Mod}~:~IM=0 \right\}$. We claim that $\mathfrak{T}$ has the closure properties of a torsion class. It is easy to see that $\mathfrak{T}$ is closed under homomorphic images and coproducts. We show that if the ideal $I$ of a ring $R$ is flat as an $R$-module, then $\mathfrak{T}$ is also closed under extensions. Let $0\rightarrow A\rightarrow B \rightarrow C\rightarrow 0$ be a short exact sequence of $R$-module with $A, C\in \mathfrak{T}$, i.e., $IA=A$ and $IC=C$. Since $I$ is a flat $R$-module, $IM\cong I \otimes_R M$ for all $R$-modules $M$ and the functor $I\otimes_R -$ is exact. Therefore, we get both the short exact sequence $0\rightarrow IA\rightarrow IB \rightarrow IC\rightarrow 0$ and the commutative diagram (\ref{cd}). In light of the flatness of $I$, the canonical monomorphisms $\alpha$ and $\beta$ are isomorphisms.
\begin{equation}\label{cd}
\begin{tikzcd}
0       \arrow[r] & IA \arrow[r, "i_A"] \arrow[d, "\cong"', "\alpha"] & IB \arrow[r, "p_B"] \arrow[d, "j"] & IC \arrow[r] \arrow[d, "\cong", "\beta"'] & 0  \\
0 \arrow[r] & A \arrow[r, "f"]                                  & B \arrow[r, "g"]                  & C \arrow[r]                               & 0
\end{tikzcd}
\end{equation}

By the Short Five-Lemma, $B\cong IB$.
Since \(j:IB\to B\) is the canonical monomorphism, the fact that \(j\) is an isomorphism implies that \(IB=B\) and $B\in \mathfrak{T}$. Let \(T\in\mathfrak{T}\) and \(F\in\mathcal{T}\). Then $IT=T$ and $IF=0$.  Let  $f:T\to F$. Then
$f(T)=f(IT)\subseteq IF=0$.  Hence
$ \operatorname{Hom}_R(T,F)=0$.
Now let $M\in\mathcal{B}_I$. Define $t(M)=IM$. Since \(M\) is \(I\)-coreduced, we have \(I^2M=IM\). Hence $It(M)=I(IM)=I^2M=IM=t(M)$, so \(t(M)\in\mathfrak T\).
Also, $I(M/t(M))=0$, hence $M/t(M)\in\mathcal{T}$. Thus every \(M\in\mathcal{B}_I\) fits into the short exact sequence $0\to IM\to M\to M/IM\to 0$ with $ IM\in\mathfrak{T}$ and $ M/IM\in\mathcal{T}$.
 Together with the fact that
$ \operatorname{Hom}_R(T,F)=0$ for all \(T\in\mathfrak T\) and \(F\in\mathcal T\), this proves that
$(\mathfrak T,\mathcal T)$
is a torsion theory on $\mathcal B_I$.
\end{enumerate}

\end{prf}

\begin{rem}\rm\noi
 The assumption that  $\mathcal{A}_I$ and $\mathcal{ B}_I$ are Serre subcategories guarantees that the functors $\text{Hom}_R(R/I, -)$ and $M\mapsto IM$ are endofunctors of  $\mathcal{A}_I$ and $\mathcal{B}_I$, respectively. Indeed, Serre subcategories are closed under subobjects, quotients and extensions. Consequently, radicals, coradicals and the associated torsion theories are defined internally in these categories, allowing the use of Stenstr\"{o}m's correspondence between left exact radicals and hereditary torsion theories.
\end{rem}


 \begin{thm}\label{crad2}
 Let $I$ be an ideal of a ring $R$. The following statements are equivalent:
  \begin{enumerate}
  \item The functor $\text{Hom}_R(R/I, -)$ is a radical on the category $R$-Mod.

  \item Every $R$-module is $I$-reduced.

  \item $I$ is an idempotent ideal.

  \item $\{M\in R~\text{-Mod}~:IM=0\}$ is a TTF.

 \item The functor  $\delta_I(M):=IM$ is an idempotent  radical on the category $R$-Mod.

  \item  Every $R$-module is $I$-coreduced.

  \end{enumerate}

\end{thm}

\begin{prf}The equivalence of 1 and 2 (resp. 5 and 6) follows from Proposition \ref{t3}(1) (resp. Proposition \ref{t3}(2)) and by taking $\mathcal{A}_I =R$-Mod (resp. $\mathcal{B}_I =R$-Mod). The equivalence of 3 and 4 is just the Jans' correspondence.  Suppose that 3 holds, i.e., $I^2=I$, then for every $R$-module $M$, $I^2M=IM$ establishing 6. Conversely, if every $R$-module is $I$-coreduced, then $R$ is also an $I$-coreduced $R$-module. So, $I^2R=IR$ which implies that $ I^2=I$. Thus, $3\Leftrightarrow 6$. Again, assume 3, i.e., $I^2=I$. Then for all $R$-modules $M$ and $m\in M$, $I^2m=0$ implies that $Im=0$ so that every $R$-module is $I$-reduced which is 2. Lastly, we prove that 2 implies 3. We always have \(I^2 \subseteq I\), it suffices to show that \(I \subseteq I^2\). Consider the \(R\)-module \(M = R / I^2\). Since \(I^2\) annihilates \(R/I^2\), we have \(I^2 M = 0\). By hypothesis, every \(R\)-module is \(I\)-reduced; in particular, \(M\) is \(I\)-reduced. Hence, for all \(m \in M\), \(I^2 m = 0\) implies \(I m = 0\). Since \(I^2 M = 0\) holds for every element of \(M\), we conclude \(I M = 0\). Now \(I M = 0\) means that for all \(x \in I\) and all \(r \in R\),
$ x \cdot (r + I^2) = 0 + I^2 \quad \text{in } R/I^2$, i.e., \(x r \in I^2\) for all \(r \in R\). In particular, taking \(r = 1\), we obtain \(x \in I^2\). Thus \(I \subseteq I^2\), which completes the proof.
\end{prf}

\begin{rem}\rm
 If $I$ is an idempotent ideal of $R$, then $\mathcal{A}_I=\mathcal{B}_I= R\text{-Mod}$ and Theorem \ref{main} retrieves Jans' correspondence given in \cite[Corollary 2.2]{Jan}.
\end{rem}

\begin{propn}\label{vasc} Let $I=(a)$ be a principal ideal of a ring $R$. The following statements hold.

\begin{enumerate}

 \item Any finitely generated $I$-coreduced $R$-module is $I$-reduced.
 \item An $I$-reduced finitely generated $R$-module is $I$-coreduced if and only if, $R$ is a zero-dimensional ring, i.e., if every prime ideal of $R$ is maximal.

\end{enumerate}

\end{propn}

\begin{prf}
 Define an $R$-endomorphism $f_a$ of the $R$-module $aM$ by $f_a(am)=a^2m$.  Note that $M$ is $I$-reduced (resp. $I$-coreduced) if and only if, $f_a$ is injective (resp. surjective). However by \cite[Lemma 3, page 23]{B}, every surjective endomorphism of a finitely generated $R$-module is an isomorphism. This proves part 1). Part 2) follows from the fact that an injective endomorphism of a finitely generated $R$-module is an isomorphism if and only if, $R$ is a zero-dimensional ring, see \cite{Vas}.
\end{prf}

\begin{exam}\rm
We obtain an  example from an idempotent-free but flat ideal over a domain.  Let $R = k[x]$ and $I = (x)$, where $k$ is a field.  Then $I$ is flat as an $R$-module but  $I^2 = (x^2) \neq (x) = I$.
Take $\mathcal{A}_I = \mathcal{B}_I = \{ M \in R\text{-Mod} \mid IM = 0 \}$. Equivalently, $\mathcal{A}_I = \mathcal{B}_I \simeq R/I\text{-Mod} \simeq k\text{-Vect}$, a Serre subcategory of $R\text{-Mod}$.
Now $\mathcal{T} = \{ M \in \mathcal{A}_I \mid (0:_M I) = M \} = \mathcal{A}_I$,
because $IM = 0$ for every $M \in \mathcal{A}_I$.
Also $\mathcal{F} = \{ M \in \mathcal{A}_I \mid (0:_M I) = 0 \} = 0$. Thus $(\mathcal{T}, \mathcal{F}) = (\mathcal{A}_I, 0)$ is a hereditary torsion theory on $\mathcal{A}_I$.
On the other hand, $\mathfrak{T} = \{ M \in \mathcal{B}_I \mid IM = M \}$. But every $M \in \mathcal{B}_I$ satisfies $IM = 0$, so $IM = M$ only when $M = 0$. Hence $\mathfrak{T} = 0$.
Therefore $(\mathfrak{T}, \mathcal{T}) = (0, \mathcal{B}_I)$ is a torsion theory on $\mathcal{B}_I$.

Theorem 3.2 thus yields $(\mathcal{A}_I, 0)$ on $\mathcal{A}_I$ and $(0, \mathcal{B}_I)$ on $\mathcal{B}_I$.
The class $\mathcal{T} = \mathcal{A}_I = \mathcal{B}_I$  occurs simultaneously as the torsion class in the first torsion theory and the torsionfree class in the second.
\end{exam}

\begin{exam}\rm
Let $k$ be a field and set
$R = k \times k \times k[x]$, $ I = k \times 0 \times (x)$. Then $I$ is flat as an $R$-module, since $I \cong k \times 0 \times (x)$, and $(x) \cong k[x]$ is a free $k[x]$-module. However,
$I^2 = k \times 0 \times (x^2) \neq k \times 0 \times (x) = I$, so $I$ is flat but not idempotent.  Now $R\text{-Mod} \cong k\text{-Vect} \times k\text{-Vect} \times k[x]\text{-Mod}$. Define $\mathcal{A}_I = \mathcal{B}_I := k\text{-Vect} \times k\text{-Vect} \times 0$.
This is a Serre subcategory of $R\text{-Mod}$.
An object of $\mathcal{A}_I = \mathcal{B}_I$ has the form $M = (U, V, 0)$,
where $U, V$ are $k$-vector spaces. The ideal $I$ acts as follows: $I(U, V, 0) = (U, 0, 0)$.
Also, $I^2(U, V, 0) = I(U, 0, 0) = (U, 0, 0)$.
Hence every object of $\mathcal{A}_I = \mathcal{B}_I$ is both $I$-reduced and $I$-coreduced.
Now compute the three classes in Theorem 3.2.
First, $\mathcal{T} = \{ M \in \mathcal{A}_I \mid IM = 0 \}$.  For $M = (U, V, 0)$,  $
IM = (U, 0, 0)$. Thus $IM = 0$ if and only if, $U = 0$. Hence $\mathcal{T} = \{ (0, V, 0) \mid V \in k\text{-Vect} \}$. Next, $\mathcal{F} = \{ M \in \mathcal{A}_I \mid (0 :_M I) = 0 \}$. For $M = (U, V, 0)$, the elements killed by $I$ are exactly the second‑component elements. Hence
$(0 :_M I) = (0, V, 0)$. Thus $(0 :_M I) = 0$ if and only if, $V = 0$. Therefore
$\mathcal{F} = \{ (U, 0, 0) \mid U \in k\text{-Vect} \}$. So the first torsion theory is
$ (\mathcal{T}, \mathcal{F}) = \bigl( \{ (0, V, 0) \}, \{ (U, 0, 0) \} \bigr)$. This is non-trivial because both classes are nonzero. For the second torsion theory, $\mathfrak{T} = \{ M \in \mathcal{B}_I \mid IM = M \}$. Since $IM = (U, 0, 0)$, we have $IM = M$ if and only if, $V = 0$. Hence $\mathfrak{T} = \{ (U, 0, 0) \mid U \in k\text{-Vect} \}$. Therefore the second torsion theory is $(\mathfrak{T}, \mathcal{T}) = \bigl( \{ (U, 0, 0) \}, \{ (0, V, 0) \} \bigr)$. Thus $\mathcal{T}$ occurs as the torsion class in the first torsion theory and as the torsionfree class in the second torsion theory: $
(\mathcal{T}, \mathcal{F}) = \bigl( \{ (0, V, 0) \}, \{ (U, 0, 0) \} \bigr)$, and $
(\mathfrak{T}, \mathcal{T}) = \bigl( \{ (U, 0, 0) \}, \{ (0, V, 0) \} \bigr)$. This illustrates Theorem 3.2 with $I$ flat but not idempotent. This example arises from a central idempotent decomposition of the ring. Consequently, the associated torsion theories are centrally split.
\end{exam}

 \section{Improvement of earlier results}

 \begin{paragraph}\noi \rm Aydogdu and Herbera in \cite[Proposition 2.7]{Herbera} give a    condition  for   the functor $\text{Hom}_R(R/I, -)$ to be a radical, namely; $R/I$ ought to be a flat $R$-module.  Given an inclusion $I\hookrightarrow R$.   $R/I$ being flat  as an $R$-module implies that applying the functor $R/I\otimes -$ gives another inclusion $I/I^2\hookrightarrow R/I$.  However, the image of $I/I^2$ in $R/I$  is zero. So, $I=I^2$ and every $R$-module is $I$-reduced.   So, we have
 $$R/I~\text{flat}~\Rightarrow I~\text{is idempotent}~\Rightarrow ~\text{there exists}~\mathcal{A}_I~ \text{consisting of}~I\text{-reduced}~R\text{-modules}$$ and these implications are in general irreversible.  It follows that existence of a non-trivial Serre  subcategory $\mathcal{A}_I$ consisting of $I$-reduced $R$-modules is a more general condition for the functor $\text{Hom}_R(R/I, -)$ to be a radical.
 \end{paragraph}

\paragraph{\bf A pro-zero inverse system and a weakly proregular ideal}
 Let ${\bf r} = (r_1, \cdots, r_n)$ be a sequence of elements of a ring $R$. To this sequence, we associate the Koszul complex $K(R; {\bf r})$.
For each $i\geq 1$,  let ${\bf r}^i$ be the sequence $(r_1^i, \cdots, r_n^i)$. There is a corresponding Koszul
complex $K(R; {\bf r}^i)$. Recall that an inverse system of $R$-modules $\{M_i\}_{i\geq 1}$ is called {\it pro-zero} if for every $i$
there is some $j\geq i$ such that the $R$-homomorphism $M_j\rightarrow M_i$ is zero.

\begin{defn}\rm
 A finite sequence ${\bf r} = (r_1, \cdots, r_n)$ in a ring $R$ is {\it weakly proregular} if, for
every $p <0$, the inverse system of $R$-modules $\{H^p(K(R; {\bf r}^i))\}_{i\geq 1}$ is  pro-zero.
\end{defn}

 \begin{defn}\rm
  An ideal is {\it weakly proregular}  if it  is generated by a weakly proregular sequence.
 \end{defn}

 \subparagraph{}
  In \cite[Proposition 4.1]{David}, it is proved that an idempotent finitely generated ideal is strongly idempotent if  and only if it is weakly proregular. We improve one of these implications by dropping the adjective; strongly.

 \begin{propn}\label{kkk}
  Any idempotent finitely generated  ideal of a ring $R$ is weakly proregular.
 \end{propn}

 \begin{prf}
  Suppose that $I$ is a finitely generated ideal of $R$ and $I^2=I$. By \cite[Example 3.6]{David}, every $R$-module $M$ is $I$-reduced. So,
  $\Gamma_I(M)= \bigcup_{k\in \Z^+}(0:_MI^k)=(0:_MI)\cong \text{Hom}_R(R/I, M)$ for all $M\in R\text{-Mod}$ and $R^q\Gamma_I(M)\cong \text{Ext}_R^q(R/I, M)$ for any $M\in R\text{-Mod}$ and for any integer $q>0$. It follows  that for any injective $R$-module $M$,   $R^q\Gamma_I(M)=0$, i.e., the functor $\Gamma_I$ is weakly stable, see \cite[Definition 2.1(1)]{Vyas}. By \cite[Theorem 0.3]{Vyas}, $I$ is weakly proregular.

 \end{prf}

\begin{propn}\label{P8}
 For any $R$-module $M$ and  an ideal  $I$ of $R$, the following statements are equivalent:

 \begin{enumerate}
  \item  $IM=0$,
  \item $(0:_MI)=M$,
  \item  $M$ is $I$-torsion and $I$-reduced,
  \item $M$ is $I$-complete and $I$-coreduced,
   \item $M\cong M/IM$.

 \end{enumerate}

\end{propn}

\begin{prf}
 $1\Leftrightarrow 2 \Leftrightarrow 5$ are trivial. Suppose $IM=0$. Then $M$ is $I$-torsion and $I$-reduced. Suppose that $M$ is $I$-torsion, i.e., $\Gamma_I(M)=M$. If $M$ is also $I$-reduced, then $I\Gamma_I(M)=0$ and $IM=0$. So, $1\Leftrightarrow 3$. Suppose that $M$ is $I$-complete and $I$-coreduced, i.e., $\Lambda_I(M)\cong M$ and $\Lambda_I(M)\cong M/IM$. Then $M\cong M/IM$. This shows that $4\Rightarrow 5$. Conversely, if $M\cong M/IM$, then  $IM=0$. So, $I^kM=IM=0$ and $\Lambda_I(M)=\lim_{\stackrel{\leftarrow}{k}} M/I^kM\cong M$  and hence $5\Rightarrow 4$.
\end{prf}

\begin{rem}\rm Proposition \ref{P8} shows that, what was proved in \cite{David} as an equivalence between a full subcategory of $R$-Mod consisting of all $R$-modules which are both $I$-reduced and $I$-torsion and a full subcategory of $R$-Mod which consists of all $R$-modules which are both $I$-coreduced and $I$-complete is actually an equality of the two subcategories.  It tells us which $I$-torsion modules are $I$-complete and vice-versa.
\end{rem}

\paragraph{} On the Serre subcategory $\mathcal{A}_I$
 of $I$-reduced modules, the torsion functor $\Gamma_I$ coincides with the radical $\text{Hom}_R(R/I,-)$. Consequently, the hereditary torsion theory of Theorem  \ref{main} can be interpreted as the torsion theory naturally associated to local cohomology with support in $I$.

\begin{exam}\rm
The assumption that $I$ is idempotent in Proposition \ref{kkk} cannot be omitted. Let $R=k[x]$ and $I=(x)$,
where $k$ is a field. Then $I^2=(x^2)\neq (x)=I$,
so $I$ is not idempotent. Since $R$ is a domain, we have $\Gamma_I(R)= \bigcup_{n\geq 1}(0:_R I^n)=0$.
Equivalently, $\operatorname{Hom}_R(R/I,R) \cong
(0:_R I) = 0$. On the other hand, the local cohomology of \(R\) with respect to \(I\) is computed by the \v{C}ech complex $0\longrightarrow R\longrightarrow R_x\longrightarrow 0$, and hence $H_I^0(R)=0$, $
H_I^1(R)\cong R_x/R$. Since $R_x/R\neq 0$, it follows that  $H_I^1(R)\neq 0$. Therefore,  $\Gamma_I(R)=0$ while $H_I^1(R)\cong R_x/R\neq 0$.  This shows that, for non-idempotent ideals, local cohomology need not collapse to degree zero. In particular, the functor
$\operatorname{Hom}_R(R/I,-)$ does not capture all of the $I$-torsion phenomena, and higher local cohomology modules may occur.
\end{exam}

\paragraph{}
 The preceding results show that idempotent ideals represent the extremal situation in which the $I$-torsion functor is already exact and local cohomology collapses to degree zero. The generalized Jans' correspondence developed in Section 3 may therefore be viewed as a categorical mechanism for recovering radical-theoretic behavior in situations where $I$ is not idempotent and higher local cohomology is non-trivial.

\section{The radical class induced by an ideal}

\paragraph{}  In this section, $S$ is a ring which is not  necessarily commutative or unital but will be regarded as an $R$-module.  For a subset $J\subseteq S$, we write  $J\vartriangleleft S$ (resp. $J\vartriangleleft_l S$) if it is an ideal of $S$   (resp. a left ideal of $S$). Let $R$ be a fixed commutative ring and $I$ an ideal of $R$. We consider the category $\mathcal{C}_R$  whose objects are rings equipped with a compatible $R$-module structure and whose morphisms are ring homomorphisms preserving the $R$-module structure. The notion of radical class below is understood relative to $\mathcal{C}_R$.

\begin{defn}\rm {\cite[Definition 2.1.1]{Wiegandt}} A class  of rings $\Psi$ is called a {\it radical class} if
\begin{enumerate}
 \item  $\Psi$ is homomorphically closed, i.e., if $S\in \Psi$ and $f: S\rightarrow T$ is a ring homomorphism, then $f(S)\in \Psi$;

 \item  for every ring $S\in \Psi$, the sum $\Psi(S):= \sum \{ J \vartriangleleft S~:~ J\in \Psi\}$ is in $\Psi$;

 \item  $\Psi(S/\Psi(S))=0$ for all rings $S\in \Psi$.
 \end{enumerate}

 The ideal $\Psi(S)$ is called the {\it $\Psi$-radical} (or just the {\it radical}) of $S$.
\end{defn}

\paragraph{\bf Examples of radical classes.}
 Examples of radical classes include:  the  class of all nil rings (the K\"{o}the nil radical class),
the class of all locally nilpotent rings (the Levitzki radical class), the class of all von Neumann regular rings, and the class of rings $R$ such that $(R, \circ)$ is a group, where the operation $\circ$ is defined by $a\circ b = a+b-ab$ for any $a, b\in R$ (the Jacobson radical class).

 \begin{thm}\label{rclass}
 Let $I$ be an idempotent ideal of a ring $R$.   The class of rings $$\Psi_I:=\{S\in \mathcal{C}_R~:~ IS=0\}$$ is a radical class.
 \end{thm}

 \begin{prf}
 \begin{enumerate}
  \item   Let $S\in \Psi_I$  and $f(S)$ be the homomorphic image of $S$. $S$ is an $R$-module and, by definition of $\Psi_I$, $IS=0$. $f(S)$ is also an $R$-module and $If(S)= f(IS)=0$. This shows that $\Psi_I$ is closed under homomorphic images.

  \item Let $\Psi_I(S):=\sum\{J\vartriangleleft S~:~ J\in \Psi_I\}$. For each $J\in \Psi_I$, $J$ is an $R$-module and $IJ=0$. This implies that their sum $\Psi_I(S)$ is also an $R$-module and $I\Psi_I(S)=0$. So, $\Psi_I(S)\in \Psi_I$.

  \item The ideal $\Psi_I(S/\Psi_I(S))$ of the ring $S/\Psi_I(S)$ is given by
  $$\sum\Bigl\{J/\Psi_I(S)\vartriangleleft S/\Psi_I(S)~:~J/\Psi_I(S)~\text{is an}~R\text{-module and}~ I(J/\Psi_I(S))=\bar{0}\Bigr\}.$$
  Moreover, $I(J/\Psi_I(S))=\bar{0}\Leftrightarrow IJ\subseteq\Psi_I(S)$. It follows that $I^2J=IJ\subseteq I\Psi_I(S)=0$ since $\Psi_I(S)\in \Psi_I$ by 2. So,
  the ideal $\Psi_I(S/\Psi_I(S))$ is contained in the ring
  $$  \left(\sum\Bigl\{J\vartriangleleft S~:~ \Psi_I(S)\subseteq J, J~\text{is an}~R\text{-module and}~ IJ=0\Bigr\}\right)/\Psi_I(S)=\bar{0}. $$
  To see that this ring vanishes, it is enough to realise that
  $$\sum\Bigl\{J\vartriangleleft S~:~ \Psi_I(S)\subseteq J, J~\text{is an}~R\text{-module and}~ IJ=0\Bigr\}\in \Psi_I$$ and $\Psi_I(S)$ is the largest ideal of $S$ contained in $\Psi_I$.
 \end{enumerate}
 \end{prf}

 \begin{cor}\label{c}Let $I$ be an idempotent ideal of a ring $R$.
  For any ring $S$ which is not necessarily commutative $$\Psi_I(S)\subseteq \Gamma_I(S).$$ We get equality when $S$ is commutative.
 \end{cor}

 \begin{prf}
  $\Psi_I(S)=\sum\{J\vartriangleleft S~:~IJ=0\}\subseteq\sum\{J\vartriangleleft_{l} S~:~IJ=0\}=\{r\in S~:~Ir=0\}=\Gamma_I(S)$. If $S$ is a commutative ring, then $J$ is a left ideal of $S$ if and only if, it is a two sided ideal of $S$. So, in this case, $\Psi_I(S)=\Gamma_I(S)$.
 \end{prf}

\begin{cor}\label{crad}If $I$ is an idempotent ideal of a ring $R$, then  $(0:_RI)$ is a radical ideal of $R$ and $$ \text{Hom}_R(R/I, R)\cong (0:_RI) = \Psi_I(R)=\Gamma_I(R).$$
\end{cor}

\begin{prf} The natural isomorphism
 $\text{Hom}_R(R/I, R)\cong (0:_RI)$  is well known. $\Psi_I(R)=\Gamma_I(R)$ is immediate from Corollary \ref{c}. Since $I$ is idempotent, $\Gamma_I(R)=\{r\in R~:~Ir=0\}$ which is nothing but the ideal $(0:_RI)$.
\end{prf}

\begin{cor}
 If $I$ is a nonzero idempotent ideal of a ring $R$ and the $R$-action on $S$ is  faithful, then the  commutative rings in the radical class $\Psi_I$ are non-unital.
\end{cor}

\begin{prf} Suppose there exists a unital commutative ring \(S\in\Psi_I\).  By Corollary \ref{c}, $\Gamma_I(S)=\Psi_I(S)$ for all $S\in T$. However, $S\in \Psi_I$ if and only if, $\Gamma_I(S)=S$ if and only if, $IS=0$ if and only if,  $I=0$ since the $R$-action on $S$ is faithful. This is a contradiction since by hypothesis $I$ is a nonzero ideal. $\Psi_I$ contains no commutative unital rings whenever $I\not=0$.
\end{prf}

\section{Application to local (co)homology}

\paragraph\noi
 Let $I$ be an ideal of a ring $R$.
The local cohomology (resp. local homology) of an $R$-module $M$ with respect to the ideal $I$ is the module $$\lim_{\stackrel{\rightarrow}{k}}\text{Ext}_R^i(R/I^k, M)~~~\text{(resp.}~\lim_{\stackrel{\leftarrow}{k}}\text{Tor}^R_i(R/I^k, M))$$ which we denote by $H_I^i(M)$ (resp. $H_i^I(M)$).  Over a von Neumann regular ring $R$, every $R$-module is both $I$-reduced and $I$-coreduced. In this section, we exploit this fact to compute spectral sequences associated with both local cohomology and local homology $R$-modules, $H_I^i(M)$ and $H^I_i(M)$ respectively.

\begin{propn} Let $R$ be a von Neumann regular ring  and $I$ an ideal of $R$. For any $p, q \in \Z^+$ and $M\in R\text{-Mod}$,
 $$ H_p^I(H_q^I(M))=
 \begin{cases}
  \Lambda_I( M), & \text{for}~p=q=0; \\

  0, & \text{otherwise}
 \end{cases}
$$ and the associated Grothendieck spectral sequence is given by

$$E_{00}^2 =H_0^I(H_0^I(M))= \Lambda_I(\Lambda_I(M))\Rightarrow \Lambda_I(M) $$ and

$$E_{pq}^2 = H_p^I(H_q^I(M))\Rightarrow 0~\text{for}~p\not=0 ~\text{or}~q\not=0.$$

\end{propn}

\begin{prf} Since $R$ is a von Neumann regular ring, every $R$-module $M$ is coreduced. So,
$$H_q^I(M)\cong\text{Tor}_q^R(R/I, M)= \begin{cases}
    R/I\otimes  M, & q=0 \\

    0, & \text{otherwise.}                                                         \end{cases}$$

    The isomorphism holds because if $M$ is $I$-coreduced, by \cite[Proposition 2.3]{David},\\ $\Lambda_I(M)$
    $\cong R/I\otimes M$. Since $R/I\otimes -$ is right exact, we have $$H_q^I(M)=L_q(\Lambda_I(M))\cong \text{Tor}_q^R(R/I, M).$$
    The   equality holds because every $R$-module over a von Neumann regular ring is flat. So, $R/I\otimes -$ is an exact functor and $\text{Tor}_i^R(R/I, M)=0$ for all $i>0$.

    $$H_p^I(H_q^I(M))=
    \begin{cases}

      H^I_p(R/I\otimes M), & q=0\\

      0, & \text{otherwise}
    \end{cases}$$

    $$ =
    \begin{cases}
     R/I \otimes (R/I \otimes M), & p=q=0\\
     0, & \text{otherwise}
     \end{cases}  ~~~~= ~~~~    \begin{cases}
     \Lambda_I(M), & p=q=0\\

     0, & \text{otherwise.}
    \end{cases}$$   The Grothendieck spectral sequence for the composition of the derived functors of \(\Lambda_I\) follows from the standard Grothendieck spectral sequence theorem.
\end{prf}

\begin{propn}\label{Hp}Let $R$ be an Artinian von Neumann regular ring  and $I$ an ideal of $R$. For any $p, q \in \Z^+$ and $M\in R\text{-Mod}$,

 $$H^p_I(H^q_I(M))=
 \begin{cases}
  \Gamma_I(M), & \text{for}~p=q=0; \\

  0, & \text{otherwise}
 \end{cases}
$$ and the associated Grothendieck spectral sequence is given by

$$E^{00}_2 =H^0_I(H^0_I(M))= \Gamma_I(\Gamma_I(M))\Rightarrow \Gamma_I(M) $$ and

$$E^{pq}_2 = H^p_I(H^q_I(M))\Rightarrow 0~\text{for}~p\not=0 ~\text{or}~q\not=0.$$

\end{propn}

\begin{prf}
 $$H^q_I(M)=\text{Ext}^q_R(R/I, M)= \begin{cases}
    \text{Hom}_R(R/I, M), & q=0 \\

    0, & \text{otherwise.}                                                         \end{cases}$$
    We prove the first equality first.  Every module of a von Neumann regular ring  $R$ is $I$-reduced. So, for all $M\in R\text{-Mod}$, $\Gamma_I(M)\cong \text{Hom}_R(R/I, M)$, \cite[Proposition 2.2(4)]{David}. Passing to the derived functors yields the first equality.     The second equality is due to the fact that an Artinian von Neumann regular ring is semisimple and  every $R$-module over a semisimple ring is projective. So, $\text{Hom}_R(R/I,  -)$ is an exact functor and $\text{Ext}_R^q(R/I,  -)=0$ for all $q>0$.

    $$H^p_I(H^q_I(M))=
    \begin{cases}

      H_I^p(\text{Hom}_R(R/I, M)), & q=0\\

      0,  & \text{otherwise}
    \end{cases} $$

    $$=
    \begin{cases}
     \text{Hom}_R(R/I, \text{Hom}_R(R/I, M)), & p=q=0\\
     0, & \text{otherwise}
     \end{cases} $$
  $$ =     \begin{cases}
     \Gamma_I(M), & p=q=0\\

     0, & \text{otherwise.}
    \end{cases}
$$   The statement on spectral sequences follows by definition.

 \end{prf}

 \begin{propn}
  Let $R$ be an Artinian von Neumann regular ring and $M$ be an $R$-module, then

  \begin{enumerate}
   \item $$H^I_p(H_I^q(M))=
   \begin{cases}
    \Gamma_I(M), & p=q=0 \cr

    0, & \text{otherwise.}
   \end{cases}$$

   \item $$H^p_I(H^I_q(M))=
   \begin{cases}
    \Lambda_I(M), & p=q=0 \cr

    0, & \text{otherwise.}
   \end{cases}$$

  \end{enumerate}

 \end{propn}

 \begin{prf}

 \begin{enumerate}
  \item

  $$H_I^q(M)=\begin{cases}
             \Gamma_I(M), & q=0 \cr
             0,   & \text{otherwise.}
            \end{cases}$$   So,
$$H^I_p(H_I^q(M))=\begin{cases}
  H_p^I(\Gamma_I(M)), & q=0 \cr
  0, & \text{otherwise}
                 \end{cases}$$

      $$= \begin{cases}
                    \Lambda_I(\Gamma_I(M)), & p=q=0 \cr
            0, & \text{otherwise}
                   \end{cases} =\begin{cases}
    \Gamma_I(M), & p=q=0 \cr
            0, & \text{otherwise.}
                   \end{cases}$$

\item Similar to that of 1) above.

\end{enumerate}
\end{prf}

Example \ref{l} below shows that, outside the von Neumann regular setting, the spectral sequence need not collapse trivially at degree zero.

\begin{exam}\label{l}\rm
Let $k$ be a field, let $R=k[x]$, $I=(x)$ and $M=R$. Then $I$ is principal and $R$ is Noetherian. The local cohomology of $R$ with respect to $I$ is computed by the \v{C}ech complex $0\longrightarrow R\longrightarrow R_x\longrightarrow 0$. Hence
$H_I^0(R)=0$ because $R$ has no nonzero $(x)$-power torsion, and
$H_I^1(R)\cong R_x/R$. Moreover, $H_I^i(R)=0$ $\text{for all } i\geq 2$.
Now $R_x/R$ is $I$-torsion, since every element of $R_x/R$ is killed by a sufficiently large power of $(x)$. Indeed, for $n\geq 1$, $x^n\left(x^{-n}+R\right)=1+R=0$. Therefore $H_I^0(R_x/R)=R_x/R$. Also, $(R_x/R)_x=0$, and hence
$H_I^1(R_x/R)=0$. Consequently, the spectral sequence
$$ E_2^{p,q}=H_I^p\bigl(H_I^q(R)\bigr) \Longrightarrow H_I^{p+q}(R)$$
has only one nonzero term, namely $E_2^{0,1}=  H_I^0\bigl(H_I^1(R)\bigr)= H_I^0(R_x/R) \cong R_x/R$. Thus $E_2^{0,1}\cong R_x/R$ survives to the limit and recovers
$H_I^1(R)\cong R_x/R$.
\end{exam}

\subsection*{Acknowledgement}
 The author was supported by the International Science Programme through the Eastern Africa Algebra Research Group. Some of the results of this paper were presented at the 2025 Conference on Rings and Polynomials at Graz University of Technology, Austria, where the author's participation was supported by Africa-UniNet P142-Uganda,
financed by the Austrian Federal Ministry of Women's Affairs, Science
and Research (BMFWF) and implemented by OeAD.  We are grateful to the Conference Organizers for all the support and to Dominic Bunnett, Alexandru Constantinescu, Martin Herschend, Dirk Kussin,  Kobi Kremnizer, Julian K\"{u}lshammer,  Bal\'{a}zs Szendr\H{o}i and Michael Wemyss for the  discussions. I am also grateful to the referee whose comments and suggestions greatly improved the paper.

\end{document}